\documentstyle[12pt]{article}

\newcommand{\sect}[1]{\setcounter{equation}{0}\section{#1}}

\def\be{\begin{equation}}
\def\ee{\end{equation}}
\def\bea{\begin{eqnarray}}
\def\eea{\end{eqnarray}}

\def\p{\varepsilon}
\def\>#1{{\bf #1}}                 
\def\1{\'{\i}}                           
\def\R{\rm I\kern-.2em R}

\begin{document}

\thispagestyle{empty}


\ 
\vspace{2cm}

\begin{center} {\LARGE{\bf{Null-plane Quantum
Universal $R$-matrix}}}

\end{center}

\bigskip\bigskip\bigskip

\begin{center}
A. Ballesteros${\dagger}$, F. J. Herranz${\dagger}$  and C. M.
Pere\~na${\ddagger}$
\end{center}

\begin{center}
\em 
{ { {}${\dagger}$ Departamento de F\1sica, Universidad de Burgos}
\\  E-09001, Burgos, Spain}

{ { {}${\ddagger}$ Departamento de F\1sica Te\'orica, Universidad
Complutense} \\   E-28040, Madrid, Spain}
\end{center}
\rm

\bigskip\bigskip\bigskip\bigskip

\begin{abstract} 
A non-linear map is applied onto the 
(non-standard) null-plane deformation of (3+1) Poincar\'e algebra
giving rise to a simpler form of this triangular quantization. A
universal
$R$-matrix for the null plane quantum algebra is then obtained from a
universal
$T$-matrix corresponding to a Hopf subalgebra. Finally, the
associated Poincar\'e Poisson--Lie group is quantized by using the
FRT approach.
\end{abstract}

\newpage


\sect {Introduction}

In a previous work a new triangular quantum deformation of the
(3+1) Poincar\'e algebra was introduced \cite{null}. This  
quantization was shown to deform in a consistent way both the
kinematical and dynamical contents \cite{Leut} of the null-plane
Poincar\'e symmetry, obtaining deformed Hamiltonians, spin and
position operators. The aim of this letter is to provide a quantum
universal $R$-matrix for this quantum Poincar\'e algebra that
completes this deformed model, including its corresponding quantum
group.

We recall that this
problem has been recently solved for the (2+1) case by
means of a contraction method and taking as starting point the
non-standard deformation of $so(2,2)$ \cite{Iran,poin}. As we
shall see in the sequel, provided a non-linear change of basis
inspired in the results given in \cite{poin} has been performed,
the (3+1) universal $R$-matrix can be obtained by using a universal
$T$-matrix technique \cite{FG,FGT}. An interesting consequence of
this procedure will be the factorized form of the universal
$R$-matrix, that is given by ordered (usual) exponentials of the
elements appearing within the corresponding classical $r$-matrix.
This kind of factorized expressions appears in a natural way   in
connection with transfer matrices problems in quantum field theory
where the $R$-matrix we obtain could be useful in order to
construct new integrable examples.

Let us perform the following transformation on the
null-plane generators
$\{{\widetilde P}_+,{\widetilde P}_i,\break {\widetilde
P}_-,{\widetilde E}_i, {\widetilde F}_i, {\widetilde
K}_3,{\widetilde J}_3\}$ $(i=1,2)$ of the quantum Poincar\'e
algebra $U_z{\cal P}(3+1)$ \cite{null}:
\bea
&&P_+={\widetilde P}_+,\qquad E_i={\widetilde E}_i,\qquad
J_3={\widetilde J}_3,\qquad i=1,2;\cr
&&P_-=e^{z{\widetilde P}_+}{\widetilde P}_-,\qquad 
P_i=e^{z{\widetilde
P}_+}{\widetilde P}_i,\cr
&&F_1=e^{z{\widetilde P}_+}({\widetilde F}_1-
z {\widetilde E}_1{\widetilde P}_- - 
 z {\widetilde J}_3{\widetilde
P}_2),\cr &&F_2=e^{z{\widetilde P}_+}({\widetilde F}_2-
z {\widetilde E}_2{\widetilde P}_- + 
z {\widetilde J}_3{\widetilde P}_1 ),\cr
&&K_3=e^{z{\widetilde P}_+}({\widetilde K}_3-
z {\widetilde E}_1{\widetilde P}_1 - 
z {\widetilde E}_2{\widetilde P}_2 ).
\label{ba}
\eea
After this change of basis, the resulting coproduct $\Delta$,
non-vanishing commutation relations, counit $\epsilon$ and antipode
$\gamma$ read 
\bea
&&\Delta(X)=1\otimes X+X\otimes 1, \qquad \mbox{for}\quad
X\in\{P_+,E_i,J_3\},\cr 
&&\Delta(Y)=1\otimes Y+Y\otimes e^{2zP_+},\quad  
\mbox{for}\quad Y\in\{P_-,P_i\},\cr 
&& \Delta(F_1)=1\otimes F_1
+F_1\otimes e^{2zP_+} - 2zP_-\otimes E_1 e^{2zP_+}
 - 2z P_2\otimes J_3 e^{2zP_+},\cr
&& \Delta(F_2)=1\otimes F_2+F_2\otimes e^{2zP_+}
 - 2zP_-\otimes E_2 e^{2zP_+} + 
2z P_1\otimes J_3 e^{2zP_+},\cr 
&& \Delta(K_3)=1\otimes K_3+K_3\otimes e^{2zP_+}
 - 2zP_1\otimes E_1 e^{2zP_+} - 
2z P_2\otimes E_2 e^{2zP_+};
\label{bb}
\eea
\bea  
&& [K_3,P_+]=\frac {e^{2zP_+} - 1}{2z},\qquad [K_3,P_-]=-P_-
-zP_1^2- zP_2^2,\cr 
&& [K_3,E_i]=E_ie^{2zP_+},\qquad [K_3,F_i]=- F_i - 2 z K_3 P_i,\cr 
&& [J_3,P_i]=-\p_{ij3}P_j, \qquad [J_3,E_i]=-\p_{ij3}E_j,
\qquad [J_3,F_i]=-\p_{ij3}F_j,\cr
&& [E_i,P_j]=\delta_{ij}\frac{e^{2zP_+} - 1}{2z}, \qquad 
[F_i,P_j]=\delta_{ij}(P_- + zP_1^2+ zP_2^2),\cr
&& [E_i,F_j]=\delta_{ij}K_3 +\p_{ij3}J_3e^{2zP_+},\qquad
[P_+,F_i]=-P_i,\cr 
&& [F_1,F_2]=2z(P_1F_2 - P_2F_1), \qquad [P_-,E_i]=-P_i;
\label{bc}
\eea 
\be
\epsilon(X) =0;\qquad  \mbox{for $X\in
\{P_\pm,P_i,E_i,F_i,K_3,J_3\}$}; 
\label{bd} 
\ee
\bea
&&\gamma(X)=-X\qquad  \mbox{for $X\in \{P_+,E_i,J_3\}$},\cr
&&\gamma(Y)=-Ye^{-2zP_+}\qquad  \mbox{for $X\in \{P_-,P_i\}$},\cr
&&\gamma(F_1)=-(F_1+2zP_- E_1 + 2 z P_2 J_3)e^{-2zP_+},\cr
&&\gamma(F_2)=-(F_2+2zP_- E_2 - 2 z P_1 J_3)e^{-2zP_+},\cr
&&\gamma(K_3)=-(K_3+2zP_1 E_1 + 2 z P_2 E_2)e^{-2zP_+}.
\label{be}
\eea

Note that both the coproduct and commutation relations are now
much simpler when compared to the original ones \cite{null}; in
particular, the quantum component $W_+^q$ of the Pauli-Lubanski
operator has no contribution in (\ref{bc}). In general, the map
(\ref{ba}) can be used to reproduce in this new basis the
physically relevant operators introduced in \cite{null}. For
instance, the deformed square of the mass $M_q^2$ is now
\be
M_q^2= P_-\frac{1-e^{-2zP_+}}z - 
P_1^2 e^{-2zP_+} - P_2^2 e^{-2zP_+},
\label{bf}
\ee
and it induces a deformed null-plane evolution governed by a 
$q$-Schr\"odinger equation that has been studied for the (2+1)
case in \cite{poin}. 

The key of our construction of the universal $R$-matrix is to
focus on the six generators appearing in the classical $r$-matrix
underlying this quantum deformation
\be
r= 2 (K_3 \wedge P_+ +E_1\wedge P_1
+ E_2\wedge P_2) ,
\label{bg}
\ee
since they close a Hopf subalgebra $U_z g$ after quantization.
The universal $T$-matrix for this Hopf subalgebra can be computed,
and this canonical element will give rise to a (factorized)
universal $R$-matrix for $U_z g$ in a straightforward way. The
important point is that this $R$-matrix can be shown to be a
universal $R$-matrix for the whole null-plane quantum Poincar\'e
algebra. As a first application of this result, the null-plane
quantum Poincar\'e group will be obtained.


\sect {Universal $T$-matrix for $U_z g$} 

The universal $T$-matrix of a Hopf algebra,
considered for the first time by Fronsdal and
Galindo \cite{FG,FGT},
is the Hopf algebra dual form
\be
T=\sum_{\mu} X^{\mu} \otimes p_{\mu},
\label{a4}
\ee
where $\{X^{\mu}\}$ is a basis for the Hopf algebra
and $\{p_{\mu}\}$ its dual:
$\langle p_{\nu},X^{\mu}\rangle =\delta_{\nu}^{\mu}$.
Let us remark that, in spite of the presence of a
particular basis in (\ref{a4}), the $T$-matrix is
by definition basis-independent.

We are interested in the Hopf algebra dual form for
$U_zg$ with generators $\{P_+,P_i,\break E_i,K_3\}$ $(i=1,2)$. Let
us choose the basis $X^{abcdef}= E_2^a E_1^b P_+^c K_3^d P_1^e
P_2^f$. Its dual basis will be given by the monomials $p_{lmnpqr}$
such that
\be
\langle p_{lmnpqr}, X^{abcdef} \rangle
=\delta_l^a\,\delta_m^b\,\delta_n^c
\,\delta_p^d\,\delta_q^e\,\delta_r^f.
\nonumber
\ee
We can express duality in an explicit way
by means of two structure tensors $E$ and $F$
that give, respectively, the product and the
coproduct in the Hopf algebra. For our purposes
it suffices to consider the latter, so we have
\bea
&&\Delta(X^{abcdef}):=
F_{ijklmn;pqrstu}^{abcdef}
\,X^{ijklmn}\otimes X^{pqrstu}
\label{a6}\\
&& p_{ijklmn}\, p_{pqrstu}=
F_{ijklmn;pqrstu}^{abcdef}\,p_{abcdef}.
\label{a7}
\eea
In order to compute the $T$-matrix  
we only need to know  very few selected components of
this tensor. This fact has been already used in \cite{NSP}, where
all the essential reasonings needed to prove the following
statements can be found. From~(\ref{a7}) and taking into account
that~$p_{000000}=1$ we  get  
\bea
&& F_{000000;pqrstu}^{abcdef}=
\delta_p^a\,\delta_q^b\,\delta_r^c\,
\delta_s^d\,\delta_t^e\,\delta_u^f\cr
&& F_{ijklmn;000000}^{abcdef}=
\delta_i^a\,\delta_j^b\,\delta_k^c\,
\delta_l^d\,\delta_m^e\,\delta_n^f\cr
&& F_{ijklmn;pqrstu}^{000000}=
\delta_i^0\,\delta_j^0\,\delta_k^0\,
\delta_l^0\,\delta_m^0\,\delta_n^0\,
\delta_p^0\,\delta_q^0\,\delta_r^0\,
\delta_s^0\,\delta_t^0\,\delta_u^0.
\eea
And by computing the coproducts of some specific
monomials $X^{abcdef}$ and comparing them to (\ref{a6}),
it can be checked that
\bea
&& F_{100000;pqrstu}^{abcdef}=
a\,\delta_{p+1}^a\,\delta_q^b\,\delta_r^c\,
\delta_s^d\,\delta_t^e\,\delta_u^f\cr
&& F_{010000;0qrstu}^{abcdef}=
b\,\delta_0^a\,\delta_{q+1}^b\,\delta_r^c\,
\delta_s^d\,\delta_t^e\,\delta_u^f\cr
&& F_{001000;00rstu}^{abcdef}=
c\,\delta_0^a\,\delta_0^b\,\delta_{r+1}^c\,
\delta_s^d\,\delta_t^e\,\delta_u^f\cr
&& F_{000lmn;000001}^{abcdef}=
f\,\delta_0^a\,\delta_0^b\,\delta_0^c\,
\delta_l^d\,\delta_m^e\,\delta_{n+1}^f\cr
&& F_{000lm0;000010}^{abcdef}=
e\,\delta_0^a\,\delta_0^b\,\delta_0^c\,
\delta_l^d\,\delta_{m+1}^e\,\delta_0^f\cr
&& F_{000l00;000100}^{abcdef}=
d\,\delta_0^a\,\delta_0^b\,\delta_0^c\,
\delta_{l+1}^d\,\delta_0^e\,\delta_0^f.
\eea
From all of them, and if we denote
\bea
&& {\hat e}_2=p_{100000}\qquad {\hat e}_1=p_{010000}
\qquad {\hat a}_+=p_{001000}\nonumber\\
&& {\hat k}_3=p_{000100} \qquad {\hat a}_1=p_{000010}
\qquad {\hat a}_2=p_{000001}
\label{a12}
\eea
 the following expression for the dual basis is derived
\be
p_{lmnpqr}=
\frac{{\hat e}_2^l}{l!}\frac{{\hat e}_1^m}{m!}
\frac{{\hat a}_+^n}{n!}\frac{{\hat k}_3^p}{p!}
\frac{{\hat a}_1^q}{q!}\frac{{\hat a}_2^r}{r!}.
\label{a13}
\ee
Then from (\ref{a4}) we get the final expression for the $U_zg$
dual form
\be
{\cal T}=e^{E_2\otimes{\hat e}_2}e^{E_1\otimes{\hat e}_1}
e^{P_+\otimes{\hat a}_+}e^{K_3\otimes{\hat k}_3}
e^{P_1\otimes{\hat a}_1}e^{P_2\otimes{\hat a}_2}.
\label{a14}
\ee


\sect {Universal $R$-matrix} 

A universal $R$-matrix for $U_z g$ can be easily deduced from the
$T$-matrix (\ref{a14}) provided there exists an algebra
homomorphism and coalgebra antihomomorphism $\Phi$ between the
quantum algebra $U_zg$ and its associated dual Hopf algebra
(quantum group)
$Fun_z(G)$ \cite{F}. If this condition is fulfilled, the element
\be
{\cal R}=(\mbox{id}\otimes \Phi){\cal T},
\label{cn}
\ee
with $\Phi$ acting on the generators of $Fun_z(G)$,  
is a solution of the quantum Yang--Baxter equation and, moreover,
\be
{\cal R}\Delta (X) {\cal R}^{-1}=\sigma\circ \Delta (X) ,
\quad \mbox{with}\quad \sigma(a\otimes b)= b\otimes a.
\label{co}
\ee

In our case, let us compute the defining relations for $Fun_z(G)$.
The elements in (\ref{a12}) are precisely dual
to the generators of~$U_zg$, so they can be taken as
generators for~$Fun_z(G)$.
The main tool in order to derive the commutation rules among them
is again the structure tensor $F$. From (\ref{a7}) we have
\be
[p_{ijklmn},p_{pqrstu}]=
(F_{ijklmn;pqrstu}^{abcdef}-
F_{pqrstu;ijklmn}^{abcdef})p_{abcdef}
\label{b4}
\ee
for two arbitrary elements in~$Fun_z(G)$.
Let us recall that we can obtain information
about the components of the tensor $F$ from the
coproduct in the quantum algebra, exactly as we did
in the previous section. In particular, the relevant terms in
order to compute
$[x,y]$, where $x$ and $y$ are two of the generators
considered above, are $X\otimes Y$ and $Y\otimes X$
($X$ and $Y$ being their respective dual elements)
appearing in the coproduct of elements of~$U_zg$~\cite{NSP}.
This, together with a careful preservation of
the order throughout the computations allow us to derive
the following relations
\be 
 [{\hat k}_3,{\hat a}_+]= 2z (e^{{\hat k}_3} -1),\qquad
 [{\hat a}_i,{\hat a}_+]= 2z {\hat a}_i \, e^{{\hat k}_3}, \qquad
 [{\hat e}_i,{\hat a}_j]= 2z \delta_{ij} (e^{{\hat k}_3} -1) ,
\label{ci}
\ee
and to show that the remaining commutators vanish.
In order to find the coproduct of $Fun_z(G)$ we consider a
$5\times 5$ matrix representation of $U_z g$:
\bea
&&D(P_+)=\frac 12 (e_{10}+ e_{40}),\quad 
D(P_1)= e_{20} , \quad  D(P_2)= e_{30} ,\quad
 D(K_3)= e_{14}+ e_{41} ,\cr
&&D(E_1)= \frac 12 (e_{12}+ e_{21}-e_{24}+ e_{42}) , \quad  
D(E_2)= \frac 12 (e_{13}+ e_{31}-e_{34}+ e_{43}),
\label{cj}
\eea
where $e_{ij}$ is the matrix with a single 1 entry at row $i$,
column $j$, and zeros at the remaining entries (note that this
representation is a classical one). By applying (\ref{cj}) in the
universal $T$-matrix (\ref{a14}) we obtain an element of the
quantum group $Fun_z(G)$:
\bea
&& D({\cal T})=\exp\{D(E_2)  {\hat e}_2\}\exp\{D(E_1)
 {\hat e}_1\}
\exp\{D(P_+) {\hat a}_+\}\cr
&&\qquad \times \exp\{D(K_3)  {\hat k}_3\}
\exp\{D(P_1) {\hat a}_1\}\exp\{D(P_2)  {\hat a}_2\} \cr
&&  
= \left(\begin{array}{ccccc}
 1& 0 & 0 & 0 & 0 \\
  \frac 12({\hat a}_+ +{\hat e}_1{\hat a}_1 +
{\hat e}_2{\hat a}_2) & 
\cosh({\hat k}_3)+ f & \frac 12 {\hat e}_1 &
\frac 12 {\hat e}_2 & 
\sinh({\hat k}_3)- f \\
 {\hat a}_1& \frac 12 {\hat e}_1 e^{-{\hat k}_3} & 1 & 0 & -
 \frac 12 {\hat e}_1 e^{-{\hat
k}_3} \\
 {\hat a}_2& \frac 12 {\hat e}_2 e^{-{\hat k}_3} & 0 & 1 & - 
\frac 12 {\hat e}_2 e^{-{\hat
k}_3} \\
  \frac 12({\hat a}_+ +{\hat e}_1{\hat a}_1 +
{\hat e}_2{\hat a}_2) & 
\sinh({\hat k}_3)+ f & \frac 12 {\hat e}_1 &
\frac 12 {\hat e}_2 & 
\cosh({\hat k}_3)- f  
 \end{array}\right)
\label{ck}
\eea
where
\be
f=\frac 18 ( {\hat e}_1^2 + {\hat e}_2^2 )\, e^{-{\hat k}_3}.
\label{cl}
\ee
Hence, the coproduct of $Fun_z(G)$ can be derived from the
multiplication of two quantum matrices
$D({\cal T})\dot\otimes D({\cal T})$:
\bea
&&\Delta({\hat k}_3)=1\otimes  {\hat k}_3 +
{\hat k}_3 \otimes 1,\quad
\Delta({\hat a}_i)=1\otimes  {\hat a}_i + 
{\hat a}_i \otimes 1,\cr
&&\Delta({\hat e}_i)= e^{{\hat k}_3}\otimes  {\hat e}_i + {\hat
e}_i \otimes 1,\cr 
&&\Delta({\hat a}_+)= e^{{\hat k}_3}\otimes  {\hat a}_+ + {\hat
a}_+ \otimes 1 - {\hat a}_1 e^{{\hat k}_3}\otimes {\hat e}_1
- {\hat a}_2 e^{{\hat k}_3}\otimes {\hat e}_2  .
\label{cm}
\eea

By recalling the coproduct and commutation rules of 
$U_zg$ given in (\ref{bb}) and (\ref{bc}), together with the
expressions (\ref{cm}) and (\ref{ci}), we get the  map we were
looking for:
\be 
 \Phi({\hat a}_+) = - 2 z K_3,\quad \Phi({\hat a}_i) = 
- 2 z E_i,\quad 
\Phi({\hat k}_3) =   2 z P_+,\quad \Phi({\hat e}_i) =  
 2 z P_i .
\label{cp}
\ee
Hence, the universal $R$-matrix for $U_zg$ is:
\bea
&&{\cal R}=\exp\{2 z E_2\otimes P_2\}\exp\{2 z E_1\otimes P_1\}
\exp\{- 2 z P_+\otimes K_3\}\cr
&&\qquad \times \exp\{2 z K_3\otimes P_+\}
\exp\{-2 z P_1\otimes E_1\}\exp\{- 2 z P_2\otimes  E_2\}. 
\label{cq}
\eea 
 
The $T$-matrix construction ensures that the element (\ref{cq}) is
a solution of the quantum Yang--Baxter equation and that (\ref{co})
is fulfilled for all the generators of the Hopf subalgebra. 
Furthermore, this condition is also true for the four remaining
generators ($P_-$, $F_1$, $F_2$ and $J_3$). This fact can be proved
by computing ${\cal R}\,\Delta \,{\cal R}^{-1}$ for each of them,
with ${\cal R}$ written in the form
${\cal R}=e^{A_1}e^{A_2}$ where
\be
A_1=2z(E_1\otimes P_1+E_2\otimes P_2  - P_+\otimes K_3),\quad
A_2=-2z(P_1\otimes E_1+P_2\otimes E_2  - K_3\otimes P_+).
\label{fa}
\ee

Therefore, we conclude that (\ref{cq}) is a universal $R$-matrix
for the null-plane deformation of Poincar\'e algebra
(\ref{bb})--(\ref{be}). This $R$-matrix is just the natural
generalization of the (2+1) result given in \cite{poin}, and can be
seen as an ordered exponentiation of the classical
$r$-matrix (\ref{bg}).


\sect {Null-plane quantum Poincar\'e group}

The following generic element of the quantum Poincar\'e group
$Fun_z(P(3+1))$:
\be
D({\cal P})=\left(\begin{array}{ccccc}
 1 & 0 & 0 & 0& 0 \\
\frac {{\hat x}^+}2+{\hat x}^- & {\hat \Lambda}_0^0 &
{\hat \Lambda}_1^0 &
{\hat
\Lambda}_2^0&
{\hat \Lambda}_3^0\\ 
{\hat x}^1 & \Lambda_0^1 & {\hat \Lambda}_1^1 &
{\hat \Lambda}_2^1&
{\hat \Lambda}_3^1 \\ 
{\hat x}^2 & {\hat \Lambda}_0^2 & {\hat \Lambda}_1^2 &
{\hat \Lambda}_2^2&
{\hat \Lambda}_3^2 \\
\frac {{\hat x}^+}2-{\hat x}^- & {\hat \Lambda}_0^3 &
{\hat \Lambda}_1^3 &
{\hat
\Lambda}_2^3 & {\hat \Lambda}_3^3
\end{array}\right),
\label{db} 
\ee
can be regarded as the non-commutative analogue of a null-plane
Poincar\'e group element obtained by means of a 
$5\times 5$ matrix representation of the quantum Poincar\'e
algebra given by the six matrices  (\ref{cj}) together with
\bea
 &&D(P_-) = e_{10} - e_{40},\qquad D(F_1)=
e_{12}+e_{21}+e_{24}-e_{42},\cr
 &&D(J_3) = e_{23} - e_{32},\qquad D(F_2)=
e_{13}+e_{31}+e_{34}-e_{43},
\label{da}
\eea
where the exponentials of the translation generators $P_\pm$, $P_i$
have been located at the left. The quantum Lorentz coordinates
${\hat \Lambda}_\nu^\mu$ satisfy   the pseudo-orthogonality
condition:  
\be
{\hat \Lambda}_\nu^\mu{\hat \Lambda}_\sigma^\rho\eta^{\nu\sigma}
=\eta^{\mu\rho},
\quad (\eta^{\mu\rho})={\mbox{diag}}\,(1,-1,-1,-1).
\label{ddc}
\ee

The coproduct of $Fun_z(P(3+1))$ is provided by $D({\cal
P})\dot\otimes D({\cal P})$. The counit and antipode come from
relations  $\epsilon (D({\cal P}))=I$ ($I$ is the $5\times 5$
identity matrix) and $\gamma(D({\cal P}))=D({\cal P})^{-1}$. On the
other hand, the associated commutation relations of the quantum
Poincar\'e group  can be deduced via the FRT approach \cite{FRT}. 
The universal $R$-matrix (\ref{cq}) written in the representation
(\ref{cj}) reduces to 
\be 
  D({\cal R})=I\otimes I +2z(D(K_3)\wedge D(P_+) 
+ D(E_1)\wedge D(P_1)+ D(E_2)\wedge D(P_2)).
\label{de}
\ee
Since this element fulfills the  property (\ref{co}) we can apply
the prescription 
\be
D({\cal R}) D({\cal P})_1 D({\cal P})_2 =
D({\cal P})_2 D({\cal P})_1 D({\cal R}),
\label{df}
\ee
where $D({\cal P})_1= D({\cal P})\otimes I$ and 
$D({\cal P})_2= I\otimes D({\cal P})$, thus obtaining the
commutation rules of $Fun_z(P(3+1))$:
\bea
&&[ {{\hat x}}^+,{{\hat x}}^i]=-2 z\, {{\hat x}}^i,\quad 
[ { {\hat x}}^+,{{\hat x}}^-]=-2 z\,{{\hat x}}^-
,\quad  [ { {\hat x}}^i,{ {\hat x}}^-]= 0,\quad  [ { {\hat x}}^1,{
{\hat x}}^2]=0,\cr  &&[ {{\hat \Lambda}}_\nu^\mu,{{\hat
\Lambda}}^\rho_\sigma]=0,\qquad
\nu,\mu,\rho,\sigma=0,1,2,3;\cr
&&[ {{\hat \Lambda}}_\nu^\mu,{ {\hat x}}^+] = -2z\,\delta_{\mu
0}({{\hat
\Lambda}}_\nu^3 -\delta_{\nu 0} + \delta_{\nu 3})-2z\,\delta_{\mu
3}({{\hat
\Lambda}}_\nu^0 +\delta_{\nu 0} - \delta_{\nu 3})\cr
&&\qquad\qquad  +z\,({\hat \Lambda}_0^\mu
+{\hat \Lambda}_3^\mu)({\hat \Lambda}_\nu^0+ {\hat
\Lambda}_\nu^3),\cr &&[ {{\hat \Lambda}}_\nu^\mu,{ {\hat
x}}^-]=\frac 12z\,\delta_{\mu 0}( -\delta_{\nu 0} +\delta_{\nu
3})+\frac 12 z\,\delta_{\mu 3}( -\delta_{\nu 0} +\delta_{\nu 3})
+\frac 12 z\,({\hat \Lambda}_0^\mu +{\hat \Lambda}_3^\mu)({\hat
\Lambda}_\nu^0-{\hat \Lambda}_\nu^3),\cr &&[ {{\hat
\Lambda}}_\nu^\mu,{ {\hat x}}^1]= z\,\delta_{\mu 2} {\hat
\Lambda}_\nu^1 +z\,\delta_{\mu 1} (-{\hat \Lambda}_\nu^0 + {\hat
\Lambda}_\nu^1 + {\hat \Lambda}_\nu^3 + \delta_{\nu 
0}-\delta_{\nu  3}) + z\,{\hat \Lambda}_\nu^1({\hat \Lambda}_0^\mu +
{\hat \Lambda}_3^\mu -1),\cr &&[ {{\hat \Lambda}}_\nu^\mu,{ {\hat
x}}^2]= z\,\delta_{\mu 1} {\hat
\Lambda}_\nu^2 +z\,\delta_{\mu 2} (-{\hat \Lambda}_\nu^0 + {\hat
\Lambda}_\nu^2 + {\hat \Lambda}_\nu^3 + \delta_{\nu 
0}-\delta_{\nu  3}) + z\,{\hat \Lambda}_\nu^2({\hat \Lambda}_0^\mu
+{\hat \Lambda}_3^\mu -1).\cr &&\label{dg}
\eea
The commutators among the quantum coordinates ${{\hat x}}^\pm$,
${{\hat x}}^i$ can be interpreted as the null-plane quantum
Poincar\'e plane. It is also worth mentioning that the
universal $R$-matrix can be used to construct a $q$-differential
calculus \cite{Castel} on the null-plane quantum Poincar\'e group;
in general, non-standard deformations exhibit interesting properties
in this context
\cite{Iranb}.

It can be checked these commutators are a Weyl quantization of the
Poisson brackets  of the coordinate functions $\{ x^\pm,x^i,
{\Lambda}_\nu^\mu\}$ on the classical Poincar\'e group, which can be
obtained by means of the Poisson bivector
\be
\{D({\cal P})\dot\otimes D({\cal P})\} =[r,D({\cal P})\dot\otimes
D({\cal P})],
\label{dh}
\ee
writing the classical $r$-matrix (\ref{bg}) in terms of the matrix
representation (\ref{cj}) and  applying the corresponding
pseudo-orthogonality relations (\ref{ddc}). This was the method
used in \cite{Mas,Lukdos} to construct the $\kappa$-Poincar\'e group
(recall that, in this case, only the (2+1) universal $R$-matrix has
been found in \cite{CGST}).


\bigskip \medskip
  
\noindent {\large{{\bf Acknowledgements}}}
 
\bigskip 
 This work has been
partially supported by DGICYT de Espa\~na  (Projects  PB92--0197 and
PB94--1115).  CMP is grateful to the
Departamento de F\1sica in Burgos for their hospitality.


\bigskip \bigskip

\footnotesize

\begin{thebibliography}{30}

\bibitem{null} A. Ballesteros, F.J. Herranz, M.A. del Olmo 
and M. Santander, {Phys. Lett. B} {351} (1995) 137.

\bibitem{Leut} H. Leutwyler and J. Stern, 
{Ann. Phys. (N.Y.)}  {112}  (1978) 94.

\bibitem{Iran}
A. Shariati, A. Aghamohammadi  and  M. Khorrami, 
 {Mod. Phys. Lett. A}  {11}  (1996) 187.

\bibitem{poin}
A. Ballesteros  and  F.J. Herranz, 
 {Mod. Phys. Lett. A}  to appear (q-alg/9605031).

\bibitem{FG} C. Fronsdal   and A. Galindo,  
{Lett. Math. Phys.} {27}  (1993) 39. 

\bibitem{FGT} C. Fronsdal   and  A. Galindo  {\it ``The Universal
$T$-Matrix"} in: Proc. of the  Joint Summer Research Conference on
conformal field theory, topological field theory and quantum groups. Holyhoke,
1992. 

\bibitem{NSP} A. Ballesteros, F.J. Herranz, M.A. del Olmo,
C.M. Pere\~na  and M. Santander,  {J. Phys. A} { 28}  (1995) 7113.


\bibitem{F}  C. Fronsdal, {\it ``Universal $T$-Matrix for
Twisted Quantum $gl(N)$"} preprint UCLA/93/TEP/3. 

\bibitem {FRT} 
 N.Y. Reshetikhin, L.A. Takhtadzhyan  and  L.D.  Faddeev,   
{Leningrad Math. J.} {1}  (1990) 193.   



\bibitem{Castel}
P. Aschieri and L. Castellani,      
  {Int. J. Mod. Phys. A}   {8}  (1993)  1667.


\bibitem{Iranb}
 A. Aghamohammadi,  M. Khorrami and A. Shariati, 
 {J. Phys. A}  {28}  (1995) L225.


 
\bibitem{Mas}
P. Ma\'slanka,   { J. Phys. A}  {26} (1993)  L1251.


\bibitem{Lukdos}
J. Lukierski  and H. Ruegg,      
  {  Phys. Lett. B}   {329}  (1994)  189.






\bibitem{CGST}
E. Celeghini, R. Giachetti, E. Sorace and M. Tarlini, 
 { J. Math. Phys.}   { 32}  (1991)  1159.




\end{thebibliography}
\end{document}